\documentclass[a4paper]{article}
\usepackage{amssymb,amscd,amsfonts,mathrsfs,amsmath}
\usepackage[all]{xy}
\input amssym.def

\def\double{\mathbb}

\def\zz{{\double Z}}

\def\rr{{\double R}}

\newtheorem{theorem}{Theorem}
\newtheorem{lemma}[theorem]{Lemma}

\newtheorem{remark}[theorem]{Remark}

\newcommand{\be}{\begin{equation}}
\newcommand{\ee}{\end{equation}}
\newcommand{\beq}{\begin{eqnarray}}
\newcommand{\eeq}{\end{eqnarray}}

\newcommand{\non}{\nonumber}

\newcommand{\cqfd}{\hfill\rule{1ex}{1ex}}

\begin{document}

\begin{center}

{\bf GROMOV'S PINCHING CONSTANT. }
\vskip 1cm
{\bf Galina Guzhvina}
\vskip 0.5cm
Mathematisches institut, Westf\"alische Wilhelms Universit\"at,\\
62 Einsteinstrasse,\\
48149 M\"unster, Germany \\[2mm]
{\tt guzhvina@math.uni-muenster.de}\\[2mm]
\today
\end{center}
\vskip 0.5cm
\begin{abstract}

\noindent  In early 80's M.Gromov showed that there exists a constant $\varepsilon$ such that any compact Riemannian manifold $M^n$ with $|K|_{M^n} \cdot diam^2(M^n) \leq \varepsilon$ can be finitely covered by a nilmanifold. The present paper illustrates by an explicit example that the pinching constant $\varepsilon$    depends on the dimension $n$ of the manifold, in particular, it decreases with the dimension at least as $\frac{12}{n^2}.$

\end{abstract}

\bigskip

\section{Introduction}

\bigskip

\noindent A compact Riemannian manifold $M^n$ is called $\varepsilon$-flat if its curvature is bounded in terms of the diameter as follows:

\bigskip

\beq	
|K| \leq \varepsilon \cdot diam^{-2}(M),\non
\eeq

\bigskip

\noindent where $K$ denotes the sectional curvature and $diam(M)$  the diameter of $M.$ If one scales an $\varepsilon-$flat metric  it remains $\varepsilon-$flat.

\bigskip

\noindent By almost flat we mean that the manifold carries $\varepsilon$-flat metrics for arbitrary $\varepsilon > 0.$ 

\bigskip

\noindent A question may arise whether there exist really substantial examples of such manifolds. In \cite{Grom} Gromov showed that any nilmanifold (= a compact quotient of a nilpotent Lie group) is almost flat. Moreover, he proved that those are up to finite quotients the only almost flat manifolds (cf. \cite{BuKa}). This is a remarkable result in the sense that it makes it possible to get information on the topology and algebraic structure of the manifold solely from assumptions on its curvature:

\bigskip

\noindent {\bf Theorem (Gromov)}

\bigskip

\noindent { \it Let $M^n$ be an $\varepsilon(n)-$ flat manifold, where 

\bigskip

\be \label{pinch}
\varepsilon(n) = \exp(-\exp(\exp  n^2)). 
\ee

\bigskip

\noindent Then $M$ is finitely covered by a nilmanifold. More precisely:

\bigskip

\noindent (i) The fundamental group $\pi_1(M)$ contains a torsion-free nilpotent normal subgroup $\psi$ of rank $n$,

\noindent (ii) The quotient $G = \pi_1(M)/\psi$ has finite order and is isomorphic to a subgroup of $O(n)$,

\noindent (iii) the finite covering of $M$ with fundamental group $\psi$ and deckgroup $G$ is diffeomorphic to a nilmanifold $N/\psi,$

\noindent (iv) The simply connected nilpotent group $N$ is uniquely determined by $\pi_1(M).$}

\bigskip

\noindent Gromov in \cite{Grom} does not specify the restrictions upon $\varepsilon(n).$
The constant (\ref{pinch}) is given as in the proof of Buser and Karcher, \cite{BuKa}. It reflects for larger $n$ what their proof can yield. 
It is clear that this constant may not be optimal: much better constants can be obtained for small $n (=2,3,4)$. In this context the problem of obtaining an effective pinching constant is, therefore, quite natural.
\noindent In particular, one may ask a question whether $\varepsilon$ should necessarily depend on the dimension of the manifold. This question we can answer in the affirmative:

\bigskip

\noindent {\bf Theorem (Main Result)}

\bigskip 

{\it \noindent In every dimension $n$ there exists a manifold $(M^n,g)$ with 

\be
|K|_{(M^n,g)} \cdot diam^2(M^n,g) < \frac{12}{n^2}
\ee

\noindent which can not be finitely covered by a nilmanifold.}

\bigskip

\section{An $n-$dimensional solvable  Lie group with the sectional curvature bounded by $\frac{12}{n^2}$. }

\bigskip

\noindent Consider a Lie group $S = \rr^n\rtimes \rr$ with the group operation $L_{(v,t)}(w,s) = (v + h(t)w,t + s)$, where $h(t) = Exp (tA), A\in GL(n,\rr).$
As can be easily seen, $S$ is solvable.

\noindent Indeed, take any two elements  $(v,t), (w,s) \in S$. Direct computation shows that their commutator is equal to $[(v,t),(w,s)] = (u, 0),$ where $u = h(-s-t)v - h(-t)v + h(-s)w - h(-s-t)w,$ and
$[[(v,t)(w,s)], (w,s)] = (h(-s)u - u, 0).$ Then $[[[(v,t)(w,s)], (w,s)], [(v,t),(w,s)]] = (0,0),$ and the conclusion follows.




\noindent Describe a lattice $\Gamma'$ in $S$.

\begin{lemma} \label{lattice}

\noindent A matrix $B \in GL(\rr,n)$ preserves a lattice in $\rr^n$ ($B\Gamma = \Gamma$) if and only if $B$ is conjugate to a matrix in $ GL(\zz,n)$.

\end{lemma}

\noindent {\it Proof}

\noindent Note first that if a matrix $B$ preserves a lattice, then any conjugate $TBT^{-1}$ of $B$ also preserves a lattice. Indeed, if $\Gamma$ is a lattice, $T\Gamma$ is also a lattice and is preserved by $TBT^{-1}$.  Let now $B$ preserves a canonical lattice  (the one spanned to an orthonormal basis of vectors). Straightforward computation shows that in this case $B$ and $B^{-1}$ have got the determinant equal to 1 or -1 and integer entrees, hence, $B$ belongs to  $GL(\zz,n)$. And the other way round: direct computation shows that a matrix $B$ from $GL(\zz,n)$ preserves a canonical lattice, therefore, its conjugate preserves a given one.
\cqfd

\bigskip 

\noindent The above Lemma shows that we can choose a lattice in $S$ in the form $\Gamma' = \Gamma \rtimes \zz$, provided that $\Gamma$ is a lattice in $\rr^n$ invariant under $Exp A$. Notice also that if $\Gamma' = \Gamma \rtimes \zz$ is a lattice, $\Gamma'_h = \frac{1}{h}\Gamma \rtimes \zz$ is also a lattice for $h \neq 0$.

\begin{lemma} \label{solv}

\noindent Suppose that $Exp A$ has eigenvalues with absolute value different from $1$. A quotient manifold of $S$ by a uniform discrete subgroup (a lattice) $\Gamma' = \Gamma \rtimes \zz$ can not be covered by a nilmanifold. 

\end{lemma}

\noindent {\it Proof }

\noindent Suppose, there exists a covering of $M'= S/\Gamma'$ by a nilmanifold $N'$.  Then, by Gromov's theorem, without loss of generality, we can regard $N'$ as $N'= N/\Psi$, where $\Psi$ is a lattice in a simply connected nilmanifold $N$ and $\Psi \subset \pi_1(M')$ as a normal subgroup of finite index $k$. Since $S$ is simply connected, $\pi_1(M') \backsimeq \Gamma'$. So, according to our assumption, the nilpotent group $\Psi$ is contained in the solvable  group $\Gamma'$ as a normal subgroup of finite index $k$. 

\noindent It means that $\Gamma'^{k} = <{g^k | g\in \Gamma'}>\subset \Psi$ is nilpotent and the element $(0,k)$ is contained in $\Gamma'^{k}.$
Moreover, $\Psi': = \Psi \cap \Gamma'^{k}$  is contained in $\Psi$ as a subgroup of finite index. Since $\Psi'$ is normal in $\Gamma'^{k}$ and
$\Gamma'^{k}$ is nilpotent we can find
normal subgroups $\Psi_1 \subset .... \subset \Psi_d=\Psi'$
such that $\Psi_i/\Psi_{i-1}$ is in the center of the group $\Gamma'^{k}/\Psi_{i-1}.$
This in turn implies that
the subgroups $\Psi_i$ are invariant under the map

\beq
c: \Psi' \rightarrow \Psi' \non
\eeq

\beq
x \rightarrow (0,k)x(0,-k) \non
\eeq

\noindent In addition we know that $c$ induces the identity
on the quotient group $\Psi_i/\Psi_{i-1}$.
Notice that $c$ is a linear map on the lattice $\Psi'$.
The above properties clearly imply that the eigenvalues of this
map are 1. On the other hand the eigenvalues are given
by the eigenvalues of $Exp(kA)$. Thus the
eigenvalues of $Exp A$ are roots of unity. A contradiction.

\cqfd

\bigskip

\noindent The next aim is to estimate the sectional curvature of the Lie group $S$ endowed with a suitable left-invariant metric.

\noindent We consider a matrix  $A$ which is given as the sum of the following two matrices:

\begin{displaymath}
\left( \begin{array}{cccccc}
   \left(\begin{array}{cc}
   		\lambda_1 & 0  \\
		 0 & \lambda_1  \\ \end{array}\right)  &  &  & 0 \\
		  &  \ldots &   \\
0 &  &  \left(\begin{array}{cc}
   		\lambda_l & 0  \\
		 0 & \lambda_l  \\ \end{array}\right) & 0 \\
0 &  &  & \begin{array}{cc}
   		\lambda_{l+1} &  \\
		 & \lambda_m \\ \end{array}
		  \\ \end{array}\right) 
\end{displaymath}

\beq
+ \non
\eeq

\begin{displaymath}
\left( \begin{array}{cccccc}
   \left(\begin{array}{cc}
   		0 & \varphi_1  \\
		-\varphi_1 & 0  \\ \end{array}\right)  &  &  & 0 \\
		  &  \ldots &   \\
0 &  &  \left(\begin{array}{cc}
   		0 & \varphi_l  \\
		-\varphi_l & 0  \\ \end{array}\right) & 0 \\
0 &  &  & \begin{array}{cc}
   		0 &  \\
		 & 0  \\ \end{array}
		  \\ \end{array} \right)
\end{displaymath}

\noindent with $\lambda_i$  real.

\noindent Then $Exp(tA) =$

\begin{displaymath}
\left( \begin{array}{cccccc}
   e^{t\lambda_1}\left(\begin{array}{cc}
   		cos(t\varphi_1) & sin(t\varphi_1)  \\
		-sin(t\varphi_1)& cos(t\varphi_1)  \\ \end{array}\right)  &  &  & 0 \\
		  &  \ldots &   \\
0 &   & e^{t\lambda_l}\left(\begin{array}{cc}
   		cos(t\varphi_l) & sin(t\varphi_l)  \\
		-sin(t\varphi_l)& cos(t\varphi_l)  \\  \end{array}\right) & 0 \\
0 &  &  & \begin{array}{cc}
   		e^{t\lambda_{l+1}} & 0 \\
		0 &  e^{t\lambda_m} \\ \end{array}
		  \\ \end{array} \right)
\end{displaymath}

\bigskip

\noindent We define a standard left-invariant metric on $S$:

\beq \label{metric}
\langle v_1,v_2 \rangle_{(w,t)} = \langle (dL_{(w,t)}^{-1})_{(w,t)}v_1,(dL_{(w,t)}^{-1})_{(w,t)}v_2\rangle_{(0,0)} = \langle\tilde {h}(-t) v_1,\tilde {h}(-t)v_2\rangle_{(0,0)}
\eeq

\noindent where

\begin{displaymath}
\tilde {h}(t) =\left( \begin{array}{cc}
   h(t)  &   0 \\
 0 &    1
		  \\ \end{array} \right),
\end{displaymath}

\noindent $<\cdot, \cdot>_{(0,0)}$  is a standard Euclidean scalar product on $s \backsimeq \rr^{n+1}$  for $s$ - the Lie algebra of $S$ and $v_1,v_2 \in s $. From the explicit expression for this metric follows the obvious

\begin{lemma}

\noindent $S$ is isometric to $\tilde S$, where $\tilde S$ is a Lie group corresponding to the matrix $A$ equal to

\begin{displaymath}
\left( \begin{array}{cccccc}
   \left(\begin{array}{cc}
   		\lambda_1 & 0  \\
		 0 & \lambda_1  \\ \end{array}\right)  &  &  & 0 \\
		  &  \ldots &   \\
0 &  &  \left(\begin{array}{cc}
   		\lambda_l & 0  \\
		 0 & \lambda_l  \\ \end{array}\right) & 0 \\
0 &  &  & \begin{array}{cc}
   		\lambda_{l+1} &  \\
		 & \lambda_m \\ \end{array}
		  \\ \end{array} \right)
\end{displaymath}

\end{lemma}

\bigskip

\begin{lemma} \label{curv}

\noindent  If one puts $\lambda_{max} = max \{|\lambda_i|, i = 1,...,m\}$, the sectional curvature of $S$ is bounded from above by  $|K|_S\leq \frac{11}{4}\lambda^2_{max}$

\end{lemma}

\noindent {\it Proof}

\noindent The curvature of a left-invariant metric on a Lie group is given by

\begin{eqnarray}
\lefteqn{<R(X,Y)Y,X> = \frac {1}{4}\|(ad_X)^*(Y) + (ad_Y)^*(X)\|^2 - <(ad_X)^*(X),(ad_Y)^*(Y)> }\nonumber\\
&&\qquad - \frac {3}{4}\|[X,Y]\|^2 - \frac {1}{2}<[[X,Y],Y],X> -\frac {1}{2}<[[Y,X],X],Y> \nonumber
\end{eqnarray}

\noindent where $X$,$Y$ are left-invariant vector fields,
$$X_{(v,t)} = (dL_{(v,t)})X_e = \tilde h(t)X_e,$$
$$Y_{(v,t)} = (dL_{(v,t)})Y_e = \tilde h(t)Y_e$$

\noindent (cf., for example, \cite{ChEb}).

\noindent To simplify the computations, for metrical estimates  we can use the group $\tilde S$.
\noindent The Lie bracket for $\tilde S$ is given by
$$[X,Y] = x_0AY' - y_0AX'$$ where
$$
X=
\left( \begin{array}{cc}
   X', & x_0 \end{array}\right)
$$
$$
Y=\left( \begin{array}{cc}
   Y', & y_0 \end{array}\right)
$$
$ X',Y' \in \rr^n, x_0,y_0 \in \rr $.

\noindent So, $$ \|ad_XY\| = \|[X,Y]\|\leq \max_i| \lambda_i| \|X\|\|Y\|$$

\noindent and the same estimation holds for the  matrix of the adjoint operator $(ad_X)^*$.

\noindent Indeed, take $(ad_X)^*Y \neq 0$ and put $Z = \frac{(ad_X)^*Y}{\|(ad_X)^*Y\|}.$ Then

\beq
\|(ad_X)^*Y\| = <(ad_X)^*Y, Z> = < Y, ad_X Z> \leq \|ad_X Z\|\|Y\| \leq \max_i| \lambda_i| \|X\|\|Y\| \non
\eeq

\noindent Hence, finally, $|K|\leq \frac{11}{4}\lambda^2_{max}$.
\cqfd

\bigskip

\noindent So we see that the sectional curvature of $S$ is controlled by the eigenvalues of the matrix $A$.

\begin{remark}

\noindent If the eigenvalues of $Exp A$ satisfy the equation
\begin{displaymath}
x^n + 1 = 0
\end{displaymath}
\noindent it corresponds to the case when $S$ is flat.
\end{remark}

\noindent Consider the equation

\be \label{charpol}
x^{2k} + 3x^{k} + 1 = 0
\ee

\begin{lemma} \label{charpoly}

\noindent The left-hand side of the equation (\ref{charpol}) is the characteristic polynomial of a matrix $T' \in  GL(\zz,n) $ if $n = 2k$.

\end{lemma}

\noindent {\it Proof}

\noindent Let $T' = $

$$
\left( \begin{array}{c}
               \left. \begin{array}{ccccccc}
                      0 &  0  &  0 & \ldots  &  & 0 & 1\\
                       -1 &  0   & 0 &\ldots  &  & 0  & 0\\
                      0 & -1 &  0 &  \ldots  &  & 0  & 0\\
                     & & &\ldots & & &  \end{array}  \right\}k\\
                    \begin{array}{ccccccc}
                   0 & 0 & 0 & \ldots  & & 0 & a \\
                    & & & \ldots  & &  & \\
                     0 & 0 & 0 &\ldots & & -1 & 0 \end{array} \end{array}\right), 
$$

\noindent where $a = (-1)^{k}\cdot3.$

\noindent Direct computation shows that this matrix is indeed in $GL(\zz,n).$ From the explicit form of $T',$ the characteristic polynomial of $T'$ is exactly the polynomial on the left-hand side of the equation (\ref{charpol}). \cqfd

\begin{remark}

\noindent If $n = 2k + 1$ we consider the polynomial

\be \label{charpol1}
(x +1)(x^{2k} + 3x^{k} + 1) = 0
\ee

\noindent and the corresponding matrix $T'' \in  GL(\zz,n)$

\begin{displaymath}
T'' =\left( \begin{array}{ccccccccc}
   \qquad\left. \begin{array}{ccccccccc}
   0 &  0  &  0 & \ldots \ldots  \ldots \ldots& &  & 0 &  0 & -1 \\
     -1 &  0   & 0 &\ldots \ldots \ldots \ldots & & & 0 & 0 & 1\\
    0 & -1 &  0 &  \ldots \ldots \ldots \ldots& &  & 0 & 0 & 0\\
       & \ldots & 
     \end{array}\right \}k+1\\ 
    \begin{array}{ccccccccc}
 0 & 0 & 0 \ldots \ldots  \ldots \ldots&  & & & 0 & 0 & a_1 \\
        0 & 0 & 0\ldots \ldots   \ldots \ldots	&  & & & 0 & 0 & a_2 \\
      & \ldots    & \\
     0 & 0 & 0\ldots \ldots  \ldots \ldots & &  & & -1  & 0 & 1\\
     0 & 0 & 0 \ldots \ldots \ldots \ldots&  & & & 0 & -1 & 0\\
\end{array}
   \end{array}\right), 
\end{displaymath}

\noindent where $a_1 = (-1)^{k}\cdot3, a_2 = (-1)^{k+1}\cdot3.$

\end{remark}

\noindent Note also that the matrix $T'(T'')$ is semisimple (cf., for example, W. Greub, \cite{Greu}), hence each of its invariant subspaces has a complement invariant subspace, therefore, $T'$ can be decomposed over the reals into $2\times 2$ blocks. In particular, $T'$ is conjugate to

\begin{displaymath}
U' = \left( \begin{array}{cccc}
   e^{\lambda_1}\left(\begin{array}{cc}
   		cos\varphi_1 & sin\varphi_1  \\
		-sin\varphi_1 & cos\varphi_1  \\ \end{array}\right)  &  & 0 \\
		  &  \ldots &   \\
0 &   & e^{\lambda_l}\left(\begin{array}{cc}
   		cos \varphi_l & sin\varphi_l  \\
		-sin\varphi_1 & cos \varphi_l  \\  \end{array}\right)
		  \\ \end{array} \right)
\end{displaymath}

\noindent for $\lambda_i = ln |r_i|,$ where  $r_i$ are the roots of the characteristic equation (\ref{charpol}), and corresponding $\varphi'$s.

\noindent Straightforward computation shows that  for any $i = 1, ..., n,$ $ln |r_i| < \frac {2}{n}$.

\noindent Hence, $\lambda_{max} < \frac {2}{n}.$

\bigskip 

\noindent Now consider the matrix A corresponding to these values of $\lambda_i$ and $\varphi_i$ and let $S$ denote the corresponding group. By construction $A$ is conjugate to a matrix in $Gl(\zz, n)$. Therefore $A$ preserves a lattice $\Gamma$ and hence $S$ contains the lattice of the form $\Gamma':=\Gamma \rtimes \zz$. Clearly, $\Gamma'_h=\frac{1}{h}\Gamma \rtimes \zz$ is lattice as well

\bigskip
 
\noindent The next step is to estimate the diameter of the quotient manifold  $S/\Gamma'_h:$

\begin{lemma} \label{diam}

\bigskip

\noindent  

\be
\lim_{h\to \infty} diam ( S/ \Gamma'_h ) = diam (\rr/\zz) = 1.
\ee

\end{lemma}

\noindent {\it Proof} {

\noindent  First note that $S \diagup_{(\frac{1}{h}\Gamma \rtimes \zz)}$ fibers over  $S \diagup_{(\rr^n\rtimes \zz)}$. The natural projection

\begin{displaymath}
pr: S \diagup_{(\frac{1}{h}\Gamma \rtimes \zz)} \rightarrow S \diagup_{(\rr^n\rtimes \zz)} = S^1
\end{displaymath}
\noindent  is a Riemannian submersion (a maximal rank surjective map, preserving the lengths of vectors orthogonal to the fiber.)
It is easy to see, that the diameter of the fiber tends to zero. Thus

\begin{displaymath}
diam M^1 \rightarrow diam S^1 = 1
\end{displaymath}

\cqfd

\bigskip

\section{Proof of the Main Result}

\bigskip

\noindent In any dimension $n$ take $(M^n,g)$ equal to $(S/\Gamma'_h, g'),$ where $S/\Gamma'_h$ is a quotient of an $n$-dimensional  solvable Lie group described in Section 2 and 
$g'$ is a left-invariant metric on $S$ as in (\ref{metric}).
\noindent From Lemma \ref{charpoly} the estimation for the maximal eigenvalue of $A$ in the definition of $S$ is $\lambda_{max} \leq \frac{2}{n}.$

\noindent Now we can use Lemma \ref{curv} to estimate the sectional curvature of $(S/\Gamma'_h, g')$:

\beq
|K|_{(S/\Gamma'_h, g')} \leq \frac{11}{n^2}.
\eeq

\noindent  From Lemma \ref{diam}  we can choose an $h$ so that 

\beq
diam(S \diagup_{(\frac{1}{h}\Gamma \rtimes \zz)}) < \sqrt{\frac{12}{11}}.
\eeq

\noindent For this $h$

\begin{displaymath}
diam^2 (S /\Gamma'_h ) \cdot |K (S /\Gamma'_h)| \leq \frac{12}{n^2}.
\end{displaymath}

\noindent  Recall that, from Lemma \ref{solv}, $S /\Gamma'_h$ can not be covered by a nilmanifold. 

\noindent  So, we can conclude that the pinching constant in the Gromov's Theorem decreases with the dimension at least as $\frac{12}{n^2}$.

\end{document}